\documentclass [12pt,a4paper,reqno]{amsart}
\usepackage{amssymb}

\def\XXint#1#2#3{{\setbox0=\hbox{$#1{#2#3}{\int}$}
     \vcenter{\hbox{$#2#3$}}\kern-.5\wd0}}

\usepackage{amsfonts,amsmath,amsthm,amscd,amssymb,latexsym}

\usepackage[T1]{fontenc}
\usepackage[utf8]{inputenc} 
\usepackage{amsmath,amssymb,amsfonts,amsthm,amscd}
\usepackage{graphicx}
\usepackage{enumitem}
\usepackage{needspace}
\usepackage{lmodern}

\newcounter{lemma}[section]

\newcounter{corollary}[section]

\newcounter{remark}[section]

\newcounter{theorem}[section]

\newcounter{proposition}[section]

\newcounter{definition}[section]

\newcounter{example}[section]

\begin{document}

\title[Weak Lipschitzness]{Weakly Lipschitz mappings in higher dimensions}
\author[A. Golberg]{Anatoly Golberg}
\address{School of Mathematical Sciences, Holon Institute of Technology,
52 Golomb St., P.O.B. 305, Holon 5810201, Israel}
\email{golberga@hit.ac.il}
\thanks{ORCID: 0009-0002-8785-7463}

\begin{abstract}
We introduce and investigate the main features of weakly Lipschitz mappings whose pointwise Lipschitz constant is controlled from above by the derivative of a set function. The main results give conditions under which such mappings belong to
$W^{1,n}_{\rm loc},$ and show that the class contains metrically
quasiconformal mappings as a proper subclass. In addition, under a Rado--Reichelderfer type oscillation assumption, we derive diameter control, image-measure bounds, sufficient conditions for Lusin's $(N)$-property, and an area formula.
\end{abstract}

\maketitle

\centerline{\emph{Dedicated to Pekka Koskela on the occasion of his 65th Birthday }}

\bigskip
{\small {\textbf {2020 Mathematics Subject Classification: }
Primary: 30C65; Secondary: 46E35, 26B05}}

\bigskip
{\small {\textbf {Key words:} weakly Lipschitz mappings, ${\rm ACL}$ mappings, Sobolev class $W^{1,n},$ set functions,  Lusin $(N)$-property, quasiconformal mappings}}

\medskip

\bigskip


\section{Introduction}\label{Intro}

In one dimension, Lipschitz continuity provides a simple and robust form of absolute continuity: Lipschitz functions are absolutely continuous, differentiable almost everywhere, and have essentially bounded derivatives. In higher dimensions the picture becomes more subtle. Nevertheless, Lipschitz mappings retain many fundamental properties. They admit McShane-type and Kirszbraun-type extensions, are differentiable almost everywhere in the sense of Rademacher, are absolutely continuous on almost all lines parallel to the coordinate axes, preserve Lebesgue null sets, and do not increase the Hausdorff dimension of sets. These features make Lipschitz mappings a central object in higher-dimensional analysis and provide a natural starting point for the weaker notions studied in this paper.

Recall that a mapping $f$ on a metric space $X$ is said to be $L$-\textit{Lipschitz} (or \textit{Lipschitz}) if there is a constant $L\ge 1$ such that
\begin{equation*}
|f(x)-f(y)|\,\le\,L|x-y|
\end{equation*}
for all $x$ and $y$ in $X;$ see, e.g., \cite{Hei01}. In Euclidean domains, every Lipschitz mapping belongs to $W^{1,\infty}.$

Note that the insertion of limsup in the previous definition, i.e.,
\begin{equation}\label{eq:limLip}
L(x,f):=\,\limsup\limits_{y\to x}\frac{|f(x)-f(y)|}{|x-y|}\,\le\,\mathcal L
\end{equation}
does not bring any essential difference. Moreover, the relationship between both the constants $L$ and $\mathcal L$ is quantitative; see, e.g., \cite{CF08}.

Recently, mappings for which only the finiteness of $L(x,f)$ is required have been studied intensively; see, e.g., \cite{AG20, AG22AAMP, AG22CMFT, KR08, MRSY09, Sal11, Sal15}. This weakening may lead to mappings which fail to be even Sobolev ones, i.e., to belong to $W^{1,1}.$ It can be illustrated by a simple example. Consider a rotation in the plane $x_3=\ldots=x_n=0,$
\begin{equation*}
f(x)\,=\,\left(x_1\cos\theta-x_2\sin\theta,x_1\sin\theta+x_2\cos\theta, x_3,\ldots,x_n\right)\,,
\end{equation*}
where $\theta=1/(x_1^2+x_2^2).$
Although this mapping is ${\rm ACL}$ after the natural continuous extension to the axis (e.g. $f(0,0,x_3,\ldots,x_n)=(0,0,x_3,\ldots,x_n)$) its pointwise Lipschitz constant behaves like $(x_1^2+x_2^2)^{-1}$ near the axis and is not locally integrable; cf. \cite{AG22AAMP}. Thus finiteness of $L(x,f)$ a.e. alone does not imply even $W^{1,1}_{\rm loc}$-regularity.
Therefore, in order to obtain $W^{1,p}_{\rm loc}$-regularity, some
appropriate upper bound for $L(x,f)$ should be assumed.
Throughout the paper, we assume that
\begin{equation}\label{eq:limset}
\limsup\limits_{y\to x}\frac{|f(x)-f(y)|}{|x-y|}\,\le\,\left[\Phi'(x)\right]^{1/n}\quad\text{a.e.}\,,
\end{equation}
where $\Phi'$ is a derivative of a set function $\Phi,$ and call such $f$ weakly Lipschitz.
By the Rademacher--Stepanov theorem \cite{Step25}, condition \eqref{eq:limset} implies that the mapping $f$ is differentiable a.e.

In this paper, we establish that continuous weakly Lipschitz mappings
satisfying a natural linewise Lusin $(N)$-condition are ${\rm ACL}$ and belong to
${\rm ACL}^n.$
We also show that metrically quasiconformal mappings form a proper subclass of weakly Lipschitz ones. Recall also that an $L$-bi-Lipschitz homeomorphism is metrically
quasiconformal with $H_f(x)\le L^2.$
We also discuss additional consequences, including Rado--Reichelderfer diameter control, image-measure bounds, Lusin-type properties, and area formulas, under
natural supplementary assumptions.


\section{Preliminary notions}\label{notions}


\subsection{Set functions}\label{subsec:Setfunc} 
Let $D$ be a domain in $\mathbb R^n.$ We consider a family $\mathcal F$ of set functions $\Phi$ defined on open
subsets $A\subset D\subset\mathbb R^n$ such that $\Phi(A)\ge 0$ and
\begin{equation}\label{eq:F-subadd}
 \sum_{k=1}^q \Phi(A_k)\le \Phi(A),
\end{equation}
whenever $A_1,\dots,A_q$ are pairwise disjoint open sets with $A_k\subset A$.
It is also assumed throughout that $\Phi$ is locally finite, i.e.
$\Phi(A)<\infty$ whenever $A\Subset D$ is open.

By \cite{RR55}, there exist the lower and upper derivatives of $\Phi\in\mathcal F$
\begin{equation*}
\underline{\Phi'}(x)\,=\,\liminf\limits_{h\to 0, d(Q)<h}\frac{\Phi(Q)}{mQ}\,,\quad
\overline{\Phi'}(x)\,=\,\limsup\limits_{h\to 0, d(Q)<h}\frac{\Phi(Q)}{mQ}\,,
\end{equation*}
where the liminf and limsup are taken over all open balls or
coordinate-parallel cubes $Q\ni x$ with $d(Q)<h.$
Here and hereafter $mQ$ and $d(Q)$ mean the $n$-dimensional Lebesgue measure and the diameter of $Q,$ respectively. It is known,
\begin{enumerate}[label=(\roman*)]
\item $\underline{\Phi'}(x)$ and $\overline{\Phi'}(x)$ are Borel functions;
\item $\underline{\Phi'}(x)=\overline{\Phi'}(x)$ a.e. in $D;$
\item for any open set $A\subset D,$ $\int_A \overline{\Phi'}(x)\, dm(x)\le \Phi(A).$
\end{enumerate}


\subsection{Examples of set functions in $\mathcal F$}\label{subsec:examples-set-functions} 

Note first that \eqref{eq:F-subadd} implies \emph{monotonicity}:
if $A_1\subset A$ are open then $\Phi(A_1)\le \Phi(A)$ (take $q=1$).
Secondly, the following \emph{closure properties} also hold. Namely, if $\Phi,\Psi\in\mathcal F$ and $c\ge 0$, then $c\Phi\in\mathcal F$ and $\Phi+\Psi\in\mathcal F$.
Moreover, if $\{\Phi_j\}\subset\mathcal F$ is an increasing sequence (pointwise on open sets),
then $\Phi(A):=\sup_j \Phi_j(A)$ also satisfies \eqref{eq:F-subadd}.
If, in addition, this $\Phi$ is locally
finite, then $\Phi\in\mathcal F.$


\subsubsection{Nonnegative measures and integral (density) functionals}
Let $\mu$ be a nonnegative locally finite Borel (Radon) measure on $D$ and define
\begin{equation*}
 \Phi(A):=\mu(A),\qquad A\subset D\ \text{open}.
\end{equation*}
Then $\Phi\in\mathcal F$ because for pairwise disjoint open $A_k\subset A$,
\begin{equation*}
 \sum_{k=1}^q \Phi(A_k)=\sum_{k=1}^q \mu(A_k)
 =\mu\!\left(\bigcup_{k=1}^q A_k\right)\le \mu(A)=\Phi(A).
\end{equation*}

A particularly important subclass is obtained by integration of a nonnegative density.
If $w\in L^1_{\mathrm{loc}}(D)$ and $w\ge 0$, then
\begin{equation}\label{eq:weighted-measure}
 \Phi(A):=\int_A w(x)\,dm(x)
\end{equation}
belongs to $\mathcal F$; moreover, at Lebesgue points of $w$ (hence a.e.) the derivative
satisfies $\Phi'(x)=w(x)$. Here and hereafter $\Phi'(x)$ denotes the common value of $\underline{\Phi'}(x)$ and $\overline{\Phi'}(x).$

\smallskip
\noindent\emph{Energy-type measures.}
Many natural choices in Analysis and Geometric Function Theory fall into \eqref{eq:weighted-measure}.
For instance,
\begin{equation*}
 \Phi(A)=\int_A |\nabla u(x)|^p\,dm(x)
 \qquad (u\in W^{1,p}_{\mathrm{loc}}(D),\ p\ge 1),
\end{equation*}
and for vector-valued mappings $f=(f_1,\ldots,f_n),$
\begin{equation*}
 \Phi(A)=\int_A \|Df(x)\|^p\,dm(x)
 \qquad (f\in W^{1,p}_{\mathrm{loc}}(D,\mathbb R^n),\ p\ge 1),
\end{equation*}
where $\|Df(x)\|$ can be taken as (for example) the operator norm or the Hilbert--Schmidt norm.
In the borderline case $p=n,$ these are the familiar $n$-energies.

If an (a.e.) Jacobian determinant $J_f$ exists and is integrable locally, then also
\begin{equation*}
 \Phi(A)=\int_A J_f(x)\,dm(x)
\end{equation*}
is in $\mathcal F$ provided $J_f\ge 0$ a.e. (e.g.\ for orientation-preserving mappings).

One may also take ``finite distortion''-type densities such as
\begin{equation*}
 \Phi(A)=\int_A K(x)\,J_f(x)\,dm(x),
 \qquad K(x)\ge 0,
\end{equation*}
whenever the product $KJ_f$ is locally integrable.

\smallskip
\noindent\emph{Singular measures.}
For fixed points $x_j\in D$ and weights $a_j\ge 0,$
\begin{equation*}
\Phi(A):=\sum_{j} a_j\,\mathbf 1_A(x_j)
\end{equation*}
defines a set function in $\mathcal F$ (the atomic measure $\sum_j a_j\delta_{x_j}$).
In this case $\Phi'(x)=0$ for a.e.\ $x$ (with respect to $m$).


\subsubsection{Image-volume under an injective mapping}
Let $f:D\to\mathbb R^n$ be a continuous injective mapping, in particular
a homeomorphism onto its image. Define
\begin{equation*}
\Phi(A):=mf(A),\qquad A\subset D\ \text{open}.
\end{equation*}
If $A_k\subset A$ are disjoint, then injectivity implies $f(A_k)$ are disjoint, hence
\begin{equation*}
\sum_{k=1}^q \Phi(A_k) =\sum_{k=1}^q mf(A_k)
=m\!\left(\bigcup_{k=1}^q f(A_k)\right)\le mf(A)=\Phi(A),
\end{equation*}
so $\Phi\in\mathcal F.$
A weighted variant is also in $\mathcal F,$
\begin{equation*}
\Phi(A):=\int_{f(A)} \rho(y)\,dm(y),
\qquad \rho\in L^1_{\mathrm{loc}}(\mathbb R^n),\ \rho\ge 0.
\end{equation*}


\subsubsection{Area with multiplicity (non-injective case)}
If $f:D\to\mathbb R^n$ is not injective, the choice $\Phi(A)=mf(A)$ may fail
\eqref{eq:F-subadd}. A standard substitute is to count multiplicity,
\begin{equation}\label{eq:areawithmult}
\Phi(A):=\int_{\mathbb R^n} N(f,A,y)\,dm(y),
\end{equation}
where $N(f,A,y)\in\{0,1,2,\dots,\infty\}$ is the \emph{multiplicity function} (the number of
points $x\in A$ with $f(x)=y$).
For pairwise disjoint $A_k\subset A$ one has
$\sum_{k=1}^q N(f,A_k,y) =  N\!\left(f,\bigcup_{k=1}^q A_k,y\right) \le N(f,A,y),$
hence
\begin{equation*}
\sum_{k=1}^q \Phi(A_k)\le \Phi(A)
\end{equation*}
in the extended nonnegative sense. Thus this functional satisfies
\eqref{eq:F-subadd}. If, in addition,
\begin{equation*}
\int_{\mathbb R^n} N(f,A,y)\,dm(y)<\infty
\qquad\text{whenever } A\Subset D \text{ is open},
\end{equation*}
then $\Phi\in\mathcal F.$ When an area formula is available, \eqref{eq:areawithmult} is closely related to $\int_A J_f(x)\,dm(x).$


\subsubsection{Packing (supremum) set functions}
Let $g:(0,\infty)\to[0,\infty)$ be any function. For open $A\subset D$ define
\begin{equation*}
\Phi_g(A):=\sup\left\{\sum_{i=1}^N g(r_i)\ :\
\{B(x_i,r_i)\}_{i=1}^N \ \text{pairwise disjoint and } B(x_i,r_i)\subset A\right\}.
\end{equation*}

Then $\Phi_g$ satisfies \eqref{eq:F-subadd} in the extended
nonnegative sense, because any disjoint family of balls chosen in each
$A_k$ is still a disjoint family in $A$ whenever the sets $A_k$
are disjoint. If, in addition, $\Phi_g(A)<\infty$ whenever $A\Subset D$ is open,
then $\Phi_g\in\mathcal F.$

A mapping-driven version is obtained by setting
\begin{equation*}
\Phi_f(A):=\sup\left\{\sum_{i=1}^N d(f(B_i))^n:\ \{B_i\}\ \text{pairwise disjoint balls in }A\right\},
\end{equation*}
which again satisfies \eqref{eq:F-subadd} in the extended nonnegative
sense by the same packing argument. If, in addition,
$\Phi_f(A)<\infty$ whenever $A\Subset D$ is open,
then $\Phi_f\in\mathcal F.$


\subsubsection{A simple non-additive example using an essential supremum}
Let $h:D\to [0,\infty]$ be measurable and define
\begin{equation*}
\Phi(A):=\bigl(\operatorname*{ess\,sup}_{x\in A} h(x)\bigr)\,mA,
\qquad A\subset D\ \text{open}.
\end{equation*}
Then for disjoint $A_k\subset A$,
\begin{equation*}
\sum_{k=1}^q \Phi(A_k)
\le \bigl(\operatorname*{ess\,sup}_{A} h\bigr)\sum_{k=1}^q mA_k
\le \bigl(\operatorname*{ess\,sup}_{A} h\bigr)mA
= \Phi(A),
\end{equation*}
hence $\Phi$ satisfies \eqref{eq:F-subadd}. If, in addition,
\begin{equation*}
\left(\operatorname*{ess\,sup}_{x\in A} h(x)\right)mA<\infty
\qquad\text{whenever } A\Subset D \text{ is open},
\end{equation*}
for instance if $h\in L^\infty_{\rm loc}(D),$ then
$\Phi\in\mathcal F.$ Typically, such $\Phi$ is not additive and
therefore is not a measure.


\subsection{Lusin's $(N)$ and $(N^{-1})$-properties}\label{subsec:Lusin-N-properties} 

Let $\Omega\subset\mathbb R^n$ be an open set and $f:\Omega\to\mathbb R^n$ be a mapping. We say that $f$ possesses the \textit{Lusin $(N)$-property} on a set $\Omega^\prime\subset\Omega$
if the implication
\begin{equation*}
mE\,=\,0\qquad\Longrightarrow\qquad mf(E)\,=0
\end{equation*}
holds for each subset $E$ of $\Omega^\prime.$

Similarly the mapping $f:\Omega\to\mathbb R^n$ satisfies the \textit{Lusin $(N^{-1})$-property}
if for each $E\subset f(\Omega)$ such that $mE=0$ we have $mf^{-1}(E)=0.$

The Lusin $(N)$- and $(N^{-1})$-properties can be treated as continuities in measure. Both have extensions to the $k$-dimensional Hausdorff measures, $k=1,\ldots,n-1.$


\subsection{Absolute continuity on lines} \label{subsec:ACL} 
Recall the notions of ${\rm ACL}$ and ${\rm ACL}^p,$ $p\ge 1;$ see, e.g., \cite{GMP17, Vai71, Vuo88}.

For $j=1,\ldots,n,$ let $R^n_j=\{x\in\mathbb R^n: x_j=0\},$ and let $T_j:\mathbb R^n\to R_j^n$ be the orthogonal projection $T_jx=x-x_je_j.$ Suppose that  $D\subset\mathbb R^n$ is an open set and $u:D\to\mathbb R$ is a continuous function. The function $u$ is called \textit{absolutely continuous on lines}, abbr. ${\rm ACL},$ if for every cube $Q$ with $\overline{Q}\subset D,$ the set $A_j\subset T_jD\subset R_j^n$ of all points $z\in T_jQ$ such that the function $t\mapsto u(z+te_j),$ $z+te_j\in Q,$ is not absolutely continuous as a function of a single variable satisfies $m_{n-1} A_j=0$ for all $j=1,\ldots,n.$ In the sequel, when $u$ is not assumed continuous, ``$u$ is ${\rm ACL}$'' means that $u$ admits a representative with the above property on almost all coordinate-parallel lines. By a \emph{representative} of $u$ we mean a function $\tilde u$ such that
$\tilde u=u$ a.e. in $D$ (i.e. the two functions differ only on a set of
Lebesgue measure zero); see, e.g., \cite[p.~148]{HK14}.

Absolutely continuous functions of a single variable on a compact interval lie between the classes of Lipschitz continuous and of differentiable a.e. functions. Thus, for absolutely continuous functions the derivatives exist a.e. and are Borel-measurable. Together with Fubini's theorem this implies that an ${\rm ACL}$ function $u: D\to\mathbb R$ has partial derivatives with respect to every variable $x_1,\ldots,x_n$ a.e. (with respect to the $n$-dimensional Lebesgue measure) in $D.$

We say that an ${\rm ACL}$ function $u:D\to\mathbb R$ is ${\rm ACL}^p,$ $p\ge 1,$ if $\partial u(x)/\partial x_j\in L^p(K),$ $j=1,\ldots,n,$ whenever $K\subset D$ is compact.
A vector-valued function is said to be ${\rm ACL}$ (${\rm ACL}^p$) if and only if each coordinate function belongs to this class.

For homeomorphisms (and, moreover, for general mappings from $L^p$) the classes $W^{1,p}$ and ${\rm ACL}^p,$ $p\ge1,$ coincide; see, e.g., \cite[Thm~A.15]{HK14}.


\subsection{Finitely Lipschitz and bi-Lipschitz mappings}\label{subsec:FinitelyL}
We say that a mapping $f:D\to \mathbb R^n$ is {\it finitely
Lipschitz} if $L(x,f)<\infty$ for all $x\in D.$ The mapping $f$ is called {\it finitely
bi-Lipschitz} if
\begin{equation}\label{eq:biLip}
0\,<\,l(x,f)\,\leq\, L(x,f)\,<\,\infty
\end{equation}
for all $x\in D$, where
\begin{equation}\label{eq:lowerl}
l(x,f):\,=\,\liminf_{y\to x}\frac{|f(x)-f(y)|}{|x-y|}\,.
\end{equation}

Note that various requirements - both weaker and stronger than \eqref{eq:biLip} - interacting with Lusin-type properties provide compelling classes of mappings (such as mappings of finite length, and those defined by metric and area distortions) whose definitions rely on a purely metric approach. For these, we refer primarily to the monograph \cite{MRSY09}; cf. \cite{AG20, AG22AAMP, AG22CMFT}.

In particular, every finitely Lipschitz mapping possesses the Lusin $(N)$-property with respect to the Hausdorff measure $\mathcal{H}^k$ for all $k=1, \dots, n$ (see \cite[Lem.~10.6]{MRSY09}; cf. \cite{AG22CMFT} for the case of Riemannian manifolds).


\subsection{Metric quasiconformality}\label{subsec:Metric-qc}
Quasiconformality can be defined in a metric space; see, e.g., \cite{HK98}.

A homeomorphism $f:X\to Y$ between metric spaces $X$ and $Y$ is said to be \textit{metrically quasiconformal} if it satisfies
\begin{equation}\label{eq:Hqc}
H_f(x)\,:=\,\limsup\limits_{r\to 0} H_f(x,r)\,\le\,H\,<\,\infty
\end{equation}
for all $x\in X$ for some $H$ independent of $x,$ where $H_f(x,r)=L(x,r,f)/l(x,r,f),$
$L(x,r,f)=\sup\{d'(f(x),f(y)): d(x,y)\le r\}$ and $l(x,r,f)=\inf\{d'(f(x),f(y)): d(x,y)\ge r\}.$

In the Euclidean setting the limsup in the above definition can be replaced by liminf (cf. \cite{HK98}). For the contemporary theory of multidimensional quasiconformal mappings we refer to the monographs \cite{GMP17, HKV20, Hei01, HK14, IM01, MRSY09, Resh89, Ric93, Vai71, Vuo88}.


\subsection{Weak Lipschitzness}\label{subsec:WeakL}
We say that a mapping $f:D\to\mathbb R^n$ is \textit{weakly Lipschitz} if there is a set function $\Phi\in\mathcal F$ defined on open subsets of $D$ such that the inequality \eqref{eq:limset} holds for a.e. $x\in D.$

Every weakly Lipschitz mapping has finite pointwise Lipschitz constant a.e. in $D$. However, the converse need not hold: finiteness of the pointwise Lipschitz constant a.e. does not by itself provide the integrability control required in \eqref{eq:limset}. The example given in the Introduction illustrates this distinction. On the other hand, every metrically quasiconformal mapping is weakly Lipschitz. This and other significant properties will be presented in the next sections.


\section{Main Results}\label{MainRes}

The method goes back to Menshoff \cite{Men37}, Reshetnyak \cite[Ch.~8]{Resh89}, and Trokhimchuk \cite{Tro64}. We also refer to our papers, where this approach was applied to a wide class of mappings, from the analytic functions in $\mathbb C$ to rather general mappings in $\mathbb R^n;$ cf. \cite{Gol92, Gol02MS, Gol02RRM, Gol06, Gol07, Gol14}.

\medskip
We start with a one-dimensional absolute continuity criterion.

\begin{lemma}\label{1D-AC} 
{\em Let $I\subset\mathbb R$ be an interval, let $h:I\to\mathbb R^n$ be
continuous, and let $g\in L^1_{\rm loc}(I),$ $g\ge 0.$ Suppose that
\begin{equation*}
\limsup_{s\to t}\frac{|h(s)-h(t)|}{|s-t|}
\le g(t)
\quad\text{for a.e. }t\in I.
\end{equation*}
Assume, in addition, that each scalar component of $h$ has Lusin's
$(N)$-property on compact subintervals of $I.$ Then
\begin{equation*}
h\in {\rm AC}_{\rm loc}(I).
\end{equation*}
Moreover, for every compact interval $[a,b]\subset I,$
\begin{equation*}
|h(b)-h(a)|\le \int_a^b g(t)\,dt.
\end{equation*}}
\end{lemma}

\begin{proof}
It is enough to consider a compact interval $[a,b]\subset I.$
By the Rademacher--Stepanov theorem in one dimension, \(h\) is differentiable
a.e. on $[a,b].$ At every point of differentiability at which the assumed
limsup estimate holds, we have
\begin{equation*}
|h'(t)|\le\limsup_{s\to t}\frac{|h(s)-h(t)|}{|s-t|}\le g(t).
\end{equation*}
Thus $h'\in L^1([a,b]).$

Since each scalar component $h_i$ is continuous, has Lusin's
$(N)$-property, and satisfies $h_i'\in L^1([a,b]),$ the Banach--Zaretsky theorem, applied componentwise,
implies that $h\in {\rm AC}([a,b]).$ Here we use the Banach--Zaretsky theorem in the following standard form:
a continuous function on a compact interval which is differentiable a.e.,
has Lusin's $(N)$-property, and has derivative in $L^1,$ is absolutely
continuous; see \cite{Heil15} (cf. \cite{AG23}). Consequently,
\begin{equation*}
h(b)-h(a)=\int_a^b h'(t)\,dt,
\end{equation*}
and hence
\begin{equation*}
|h(b)-h(a)|\le \int_a^b |h'(t)|\,dt
\le \int_a^b g(t)\,dt.
\end{equation*}
Since $[a,b]\subset I$ was arbitrary, $h\in {\rm AC}_{\rm loc}(I).$
\end{proof}

\medskip
\begin{theorem}\label{thm:ACLp} 
{\em Let $f:D\to\mathbb R^n$ be a continuous weakly Lipschitz mapping whose scalar coordinate functions, restricted to almost every coordinate-parallel
line, have Lusin's $(N)$-property on compact subintervals. Then $f$ is ${\rm ACL}.$ Moreover $f\in {\rm ACL}^n.$}
\end{theorem}

\begin{proof}
Put
\begin{equation*}
g(x):=L(x,f).
\end{equation*}
By weak Lipschitzness,
\begin{equation*}
g(x)^n\le \Phi'(x)
\quad\text{for a.e. }x\in D.
\end{equation*}
Let $K\Subset D,$ and choose an open set $A$ such that
$K\Subset A\Subset D.$ By the differentiation theorem for set functions,
\begin{equation*}
\int_K g(x)^n\,dm(x)\le
\int_A \Phi'(x)\,dm(x)\le
\Phi(A)<\infty.
\end{equation*}
Thus $g\in L^n_{\rm loc}(D).$

Fix a coordinate direction $e_j$ and a parallelepiped
$P=P_0\times I\Subset D,$ whose sides are parallel to the coordinate axes.
By the Fubini theorem, for a.e. $z\in P_0,$ the function
\begin{equation*}
t\mapsto g(z+te_j)
\end{equation*}
belongs to $L^n(I),$ hence to $L^1(I).$ By the assumed linewise
Lusin $(N)$-property, for a.e. such $z,$ each scalar component of the restriction
\begin{equation*}
h_z(t):=f(z+te_j)
\end{equation*}
has Lusin's $(N)$-property on compact subintervals of $I.$

For a.e. $t\in I,$
\begin{equation*}
\limsup_{s\to t} \frac{|h_z(s)-h_z(t)|}{|s-t|}
\le L(z+te_j,f) = g(z+te_j).
\end{equation*}
By Lemma~\ref{1D-AC}, $h_z$ is absolutely continuous on $I.$
Hence $f$ is ${\rm ACL}.$

Moreover, for a.e. $x\in D,$
\begin{equation*}
\left|\frac{\partial f}{\partial x_j}(x)\right|
\le L(x,f)=g(x).
\end{equation*}
Therefore, for every $K\Subset D,$
\begin{equation*}
\int_K \left|\frac{\partial f}{\partial x_j}(x)\right|^n\,dm(x)
\le \int_K g(x)^n\,dm(x)<\infty.
\end{equation*}
Since this holds for each $j=1,\ldots,n,$ we obtain $f\in {\rm ACL}^n.$
\end{proof}

\begin{remark}
The additional linewise Lusin $(N)$-condition cannot be omitted. Indeed,
let $C$ be the Cantor function and define
\begin{equation*}
f(x_1,\ldots,x_n)=(C(x_1),x_2,\ldots,x_n)
\end{equation*}
in $(0,1)^n.$ Then $f$ is continuous and satisfies
\begin{equation*}
L(x,f)<\infty
\quad\text{for a.e. }x,
\end{equation*}
indeed $L(x,f)\le 1$ for a.e. $x.$ Hence $f$ is weakly Lipschitz with
respect to a suitable absolutely continuous set function. However, $f$ is
not ${\rm ACL},$ because its restriction to every line parallel to the $x_1$-axis
contains the Cantor function.
\end{remark}

\medskip
\begin{theorem}\label{thm:Metrqc} 
{\em Every metrically quasiconformal homeomorphism $f:D\to D'$ between domains in $\mathbb R^n$ is weakly Lipschitz with $\Phi(A)=H^n mf(A).$}
\end{theorem}

\begin{proof}
The simple chain shows for $r=|x-y|,$
\begin{equation*}
\frac{|f(x)-f(y)|}{|x-y|}\,\le\,\frac{L(x,r,f)\bigl[\Omega_n l(x,r,f)^n\bigr]^{1/n}}{l(x,r,f)\bigl[\Omega_n r^n\bigr]^{1/n}}
\,\le\, \frac{L(x,r,f)}{l(x,r,f)}\left[\frac{mf(B(x,r))}{mB(x,r)}\right]^{1/n}\,,
\end{equation*}
where $\Omega_n$ denotes the volume of the unit ball in $\mathbb R^n.$
Therefore, passing to the $\limsup$ as $r\to 0,$ we obtain
\begin{equation*}
\limsup\limits_{y\to x}\frac{|f(x)-f(y)|}{|x-y|}\,\le\,\left[\Phi'(x)\right]^{1/n}\qquad \text{a.e.}
\end{equation*}
Thus, \eqref{eq:limset} holds for $\Phi(A)=H^n mf(A),$ where $H$ is the upper bound in \eqref{eq:Hqc}.
\end{proof}

\medskip
Theorem~\ref{thm:Metrqc} shows that metrically quasiconformal mappings form a subclass of weakly Lipschitz
ones. We now give a simple example demonstrating that this inclusion is proper.

Let $D=(-1,1)^n$ and fix a parameter $\alpha\in(0,1/n)$. Define an odd, increasing
homeomorphism $\varphi:(-1,1)\to(-1,1)$ by
\begin{equation*}
\varphi(t):=\frac{1-\alpha}{2-\alpha}\left(t+\frac{\operatorname{sgn}(t)\,|t|^{1-\alpha}}{1-\alpha}\right),
\qquad t\in(-1,1),
\end{equation*}
where $\operatorname{sgn}(0)=0$. The formula extends continuously to $[-1,1],$ and this extension satisfies
\begin{equation*}
\varphi(1)=1,\qquad \varphi(-1)=-1,
\end{equation*}
hence $\varphi$ is a homeomorphism of $(-1,1)$ onto itself. Moreover, for $t\neq 0,$
\begin{equation}\label{eq:phi-derivative}
\varphi'(t)=\frac{1-\alpha}{2-\alpha}\bigl(1+|t|^{-\alpha}\bigr)>0,
\end{equation}
so $\varphi$ is strictly increasing and $\sup_{t\in(-1,1)}\varphi'(t)=\infty.$

In addition, since $\alpha n<1$ we have
\begin{equation*}
\int_{-1}^{1}\max\{1,\varphi'(t)\}^n\,dt<\infty,
\end{equation*}
because near $t=0$ the integrand is comparable to $|t|^{-\alpha n}$ and
$\int_0^1 t^{-\alpha n}\,dt<\infty.$

Now define the mapping $f:D\to D$ by the one-direction stretching
\begin{equation*}
f(x_1,\dots,x_{n-1},x_n)=(x_1,\dots,x_{n-1},\varphi(x_n)).
\end{equation*}
At points where $\varphi$ is differentiable (hence for a.e.\ $x\in D$) we have
\begin{equation*}
Df(x)=\mathrm{diag}(1,\dots,1,\varphi'(x_n)),
\end{equation*}
and, therefore,
\begin{equation*}
L(x,f)=\max\{1,\varphi'(x_n)\}\quad\text{for a.e. }x\in D.
\end{equation*}
Let
\begin{equation*}
w(x):=\max\{1,\varphi'(x_n)\}^n\in L^1_{\mathrm{loc}}(D),
\qquad
\Phi(A):=\int_A w(x)\,dm(x),\qquad A\subset D\ \text{open}.
\end{equation*}
Then $\Phi\in\mathcal F$ and $\Phi'(x)=w(x)$ for a.e.\ $x,$ whence
\begin{equation*}
L(x,f)\le [\Phi'(x)]^{1/n}\quad\text{for a.e. }x\in D,
\end{equation*}
so $f$ is weakly Lipschitz.

On the other hand, $f$ is not metrically quasiconformal, since the pointwise linear
dilatation satisfies (a.e.)
\begin{equation*}
H_f(x)\ge
\frac{\max\{1,\varphi'(x_n)\}}{\min\{1,\varphi'(x_n)\}}
=\max\{\varphi'(x_n),\,1/\varphi'(x_n)\},
\end{equation*}
which is unbounded by \eqref{eq:phi-derivative}.
Thus the class of weakly Lipschitz homeomorphisms strictly contains the class of
metrically quasiconformal mappings.


\section{Further results and remarks}\label{Otherres}

The Menshoff--Reshetnyak method together with Rado--Reichelderfer set functions typically allows one
to pursue additional conclusions of the following types (often under mild extra assumptions).


\subsection{Rado--Reichelderfer oscillation and diameter estimates}\label{subsec:RRdiam} 

In this subsection we assume, in addition, that $f$ satisfies a
Rado--Reichelderfer type estimate: there exists $\Psi\in\mathcal F$
such that
\begin{equation}\label{eq:RR-osc}
|f(x)-f(y)|^n\le \Psi(B(x,|x-y|))
\end{equation}
whenever $B(x,|x-y|)\Subset D.$

The above estimate is stronger than weak Lipschitzness. Indeed, if \eqref{eq:RR-osc} holds
for a set function $\Psi\in\mathcal F,$ then, for $r=|x-y|,$
\begin{equation*}
\left(\frac{|f(x)-f(y)|}{|x-y|}\right)^n
\le \frac{\Psi(B(x,r))}{r^n}
= \Omega_n\,\frac{\Psi(B(x,r))}{mB(x,r)}.
\end{equation*}
Passing to the limsup as $y\to x$ gives
\begin{equation*}
L(x,f)^n\le \Omega_n\Psi'(x)
\quad\text{for a.e. }x\in D.
\end{equation*}
Thus $f$ is weakly Lipschitz with respect to $\Phi=\Omega_n\Psi.$

A convenient next step is to convert \eqref{eq:RR-osc} into a \emph{uniform diameter bound}
on small balls (and hence on small cubes).

\begin{theorem}\label{thm:RRdiam} 
{\em Let $f:D\to\mathbb R^n$ be a mapping satisfying the Rado--Reichelderfer
estimate \eqref{eq:RR-osc} with some $\Psi\in\mathcal F.$ Then there exist constants
$C_1$ and $C_2$ depending on $n$ only such that the following holds:

\smallskip
\noindent\textup{(i)} For every ball $B=B(x,r)$ with $\overline{2B}\subset D$,
\begin{equation}\label{eq:diam-ball}
\bigl(\operatorname{diam} f(B)\bigr)^n \le C_1\,\Psi(2B).
\end{equation}

\smallskip
\noindent\textup{(ii)} For every cube $Q\Subset D$ with $\overline{2Q}\subset D,$
\begin{equation*}
\bigl(\operatorname{diam} f(Q)\bigr)^n \le C_2\,\Psi(2Q).
\end{equation*}
}
\end{theorem}

\begin{proof}
(i) Fix $B=B(x,r)$ with $\overline{2B}\subset D$, and take arbitrary $u,v\in B$.
Due to (\ref{eq:RR-osc}),
\begin{equation*}
|f(u)-f(x)|\le \Psi(B(x,|u-x|))^{1/n}\le \Psi(B)^{1/n},
\end{equation*}
and similarly
\begin{equation*}
|f(v)-f(x)|\le \Psi(B)^{1/n}.
\end{equation*}
Hence by the triangle inequality
\begin{equation*}
|f(u)-f(v)|\le 2 \Psi(B)^{1/n}.
\end{equation*}
Taking the supremum over all $u,v\in B$ yields ${\rm diam}\,f(B)\le 2 \Psi(B)^{1/n}.$
And by monotonicity of $\Psi$ this also implies,
\begin{equation*}
\operatorname{diam}f(B)\le 2 \Psi(2B)^{1/n}.
\end{equation*}
Therefore, \eqref{eq:diam-ball} holds with $C_1=2^n$.

\medskip
(ii) Let $Q=Q(x_Q,\ell)$ be a cube of side length $\ell,$ and fix $u,v\in Q.$
If $u=v,$ there is nothing to prove, so assume $u\neq v.$ Set
\begin{equation*}
\nu:=\left\lceil \frac{4|u-v|}{\ell}\right\rceil
\end{equation*}
and define points $x_j\in [u,v]\subset Q$ by
\begin{equation*}
x_j:=u+\frac{j}{\nu}(v-u), \qquad j=0,1,\dots,\nu.
\end{equation*}
For each $j,$ let
\begin{equation*}
B_j:=B\!\left(x_j,\frac{\ell}{4}\right).
\end{equation*}
Then $x_j\in Q,$ and for every $j,$
\begin{equation*}
2B_j=B\!\left(x_j,\frac{\ell}{2}\right)\subset 2Q\subset D.
\end{equation*}
Indeed, if $z\in B(x_j,\ell/2),$ then for each coordinate $k=1,\dots,n,$
\begin{equation*}
|z_k-(x_Q)_k| \le |z-x_j|+|(x_j)_k-(x_Q)_k|
< \frac{\ell}{2}+\frac{\ell}{2} =\ell,
\end{equation*}
hence $z\in 2Q.$

Moreover,
\begin{equation*}
|x_j-x_{j-1}|=\frac{|u-v|}{\nu}\le \frac{\ell}{4},
\end{equation*}
so the consecutive balls $B_{j-1}$ and $B_j$ intersect. Choose points
\begin{equation*}
z_j\in B_{j-1}\cap B_j,\qquad j=1,\dots,\nu,
\end{equation*}
and set $z_0=u,$ $z_{\nu+1}=v.$ Then $z_j,z_{j+1}\in B_j$ for $j=0,\dots,\nu.$ Therefore,
by part {\rm(i)},
\begin{equation*}
|f(u)-f(v)| \le \sum_{j=0}^{\nu} |f(z_j)-f(z_{j+1})|
\le \sum_{j=0}^{\nu} \operatorname{diam} f(B_j)
\le 2(\nu+1)\Psi(2Q)^{1/n}.
\end{equation*}

Since $|u-v|\le \operatorname{diam} Q=\sqrt n\,\ell,$ we have
\begin{equation*}
\nu \le \frac{4|u-v|}{\ell}+1 \le 4\sqrt n+1,
\end{equation*}
and thus
\begin{equation*}
|f(u)-f(v)| \le 2(4\sqrt n+2)\Psi(2Q)^{1/n}.
\end{equation*}
Taking the supremum over all $u,v\in Q$ and raising to the power $n,$ we obtain
\begin{equation*}
\bigl(\operatorname{diam} f(Q)\bigr)^n \le C_2\,\Psi(2Q),
\end{equation*}
where $C_2=2^n(4\sqrt n+2)^n.$ This proves {\rm(ii)}.
\end{proof}

The estimate \eqref{eq:diam-ball} is a typical Rado--Reichelderfer type bound; it can be interpreted in
terms of the Rado--Reichelderfer ${\rm RR}^n$ and Mal\'y ${\rm AC}^n$ classes; see \cite{Csor00, Hencl02, Maly99, RR55}.
In particular, bounds of the form \eqref{eq:diam-ball} control oscillations of $f$ on small balls/cubes and,
via standard covering arguments, lead to quantitative estimates for the $n$-dimensional measure of
images. This will be used in the next subsection to obtain an upper bound for $m f(A)$ in terms of
the set function $\Psi.$


\subsection{Lusin's $(N)$-property and quantitative image-measure bounds}\label{subsec:Nproperty} 

The diameter control from Theorem~\ref{thm:RRdiam} is a convenient tool for estimating the $n$-dimensional
measure of images. We first record a quantitative bound in terms of the Rado--Reichelderfer
set function $\Psi.$

\begin{proposition}\label{prop:image-measure} 
{\em Let $f:D\to\mathbb R^n$ be continuous. Assume that there exists
$\Psi\in\mathcal F$ such that the diameter estimate \eqref{eq:diam-ball} holds.
Then there exists a constant $C=C(n)\ge 1$ such that for every open set
$A \Subset D,$
\begin{equation}\label{eq:image-measure-Psi}
m f(A) \le C\,\Psi(A).
\end{equation}}
\end{proposition}

\begin{proof}
Fix an open set $A \Subset D.$ Consider the family of balls
\begin{equation*}
\mathcal B :=
\left\{B\!\left(x,\frac{1}{20}{\rm dist}(x,\partial A)\right) : x\in A
\right\}.
\end{equation*}
By the $5r$-covering lemma, there exists a countable pairwise disjoint subfamily
$\{B_i\}=\{B(x_i,r_i)\}\subset \mathcal B$ such that
\begin{equation*}
A \subset \bigcup_i 5B_i .
\end{equation*}
Since
\begin{equation*}
r_i=\frac{1}{20}{\rm dist}(x_i,\partial A),
\end{equation*}
we have
\begin{equation*}
10B_i \subset A \subset D
\qquad\text{for all }i.
\end{equation*}
Therefore, by countable subadditivity of Lebesgue measure and by Theorem~\ref{thm:RRdiam}(i)
applied to the ball $5B_i,$
\begin{equation*}
m f(A) \le \sum_i m f(5B_i)
\le \Omega_n \sum_i \bigl({\rm diam}\, f(5B_i)\bigr)^n
\le \Omega_n C_0 \sum_i \Psi(10B_i),
\end{equation*}
where $C_0=C_0(n)$ is the constant in \eqref{eq:diam-ball}.

It remains to estimate $\sum_i \Psi(10B_i).$ Suppose that
\begin{equation*}
10B_i \cap 10B_j \neq \varnothing .
\end{equation*}
Then
\begin{equation*}
20r_i={\rm dist}(x_i,\partial A)
\le |x_i-x_j|+{\rm dist}(x_j,\partial A)
< 10r_i+10r_j+20r_j,
\end{equation*}
and hence $r_i<3r_j.$ Interchanging $i$ and $j,$ we also get $r_j<3r_i.$ Thus,
whenever $10B_i$ and $10B_j$ intersect, their radii are comparable:
\begin{equation*}
\frac13 r_i < r_j < 3r_i.
\end{equation*}
Moreover, if $10B_i\cap 10B_j\neq\varnothing,$ then
\begin{equation*}
|x_i-x_j| < 10r_i+10r_j < 40r_i,
\end{equation*}
so every such center $x_j$ lies in $B(x_i,40r_i).$ Since the original balls $B_j$
are pairwise disjoint and each of them has radius $r_j \ge r_i/3,$ a standard packing
argument shows that each ball $10B_i$ meets at most $N=N(n)$ other balls $10B_j.$

A standard finite-coloring argument shows that
the family $\{10B_i\}$ can be
partitioned into $N_0=N_0(n)$ subfamilies
\begin{equation*}
\mathcal B_1,\dots,\mathcal B_{N_0},
\end{equation*}
such that the balls within each $\mathcal B_k$ are pairwise disjoint. 
Since every ball in
every $\mathcal B_k$ is contained in $A,$ property \eqref{eq:F-subadd} of the class $\mathcal F$ implies
\begin{equation*}
\sum_{10B_i\in \mathcal B_k}\Psi(10B_i)\le \Psi(A),
\end{equation*}
where the countable sum is understood as the supremum of finite partial sums. Summing
over $k=1,\dots,N_0,$ we obtain
\begin{equation*}
\sum_i \Psi(10B_i)\le N_0\,\Psi(A).
\end{equation*}
Therefore,
\begin{equation*}
m f(A)\le \Omega_n C_0 N_0\,\Psi(A)=:C\,\Psi(A),
\end{equation*}
which proves \eqref{eq:image-measure-Psi}.
\end{proof}

\medskip

A simple sufficient condition for Lusin's $(N)$-property is that $\Psi$ is \emph{absolutely continuous}
with respect to Lebesgue measure in the following outer form: null sets can be covered by open sets
of arbitrarily small $\Psi$--mass (this holds, for example, when $\Psi(A)=\int_A w\,dm$ with
$w\in L^1_{\rm loc}$, $w\ge 0$).

\medskip
\begin{corollary}\label{cor:Nproperty} 
{\em Assume the hypotheses of Proposition~\ref{prop:image-measure}. Suppose that for every measurable
set $E\Subset D$ with $mE=0$ and every $\eta>0$ there exists an open set $A\supset E$ such that
\begin{equation}\label{eq:psi-outer-ac}
\Psi(A)<\eta .
\end{equation}
Then $f$ has Lusin's $(N)$-property in $D$, i.e.
\begin{equation*}
mE=0\quad\Longrightarrow\quad m f(E) =0 .
\end{equation*}}
\end{corollary}

\begin{proof}
Let $E\subset D$ with $m E=0.$ Choose an exhaustion
$D_1\Subset D_2\Subset \cdots \Subset D$ with
\begin{equation*}
D=\bigcup_{j=1}^{\infty}D_j.
\end{equation*}
Put $E_j=E\cap D_j.$ Then $E_j\Subset D$ and $mE_j=0.$ By the assumption, for every $\eta>0$ there exists an open set
$A_j\supset E_j$ such that $\Psi(A_j)<\eta.$ Since $E_j\Subset D,$ choose an open set $U_j\Subset D$ with
$E_j\subset U_j.$ Replacing $A_j$ by $A_j\cap U_j,$ and using
monotonicity of $\Psi,$ we may assume that $A_j\Subset D.$ Hence, by
Proposition~\ref{prop:image-measure},
\begin{equation*}
m f(E_j)\le m f(A_j)\le C\Psi(A_j)<C\eta.
\end{equation*}
Letting $\eta\to 0,$ we obtain $m f(E_j)=0$ for every $j.$ Since
\begin{equation*}
E=\bigcup_{j=1}^{\infty}E_j,
\end{equation*}
we get
\begin{equation*}
m f(E)\le \sum_{j=1}^{\infty}m f(E_j)=0.
\end{equation*}
Thus $f$ has Lusin's $(N)$-property in $D.$
\end{proof}

\begin{remark}\label{rem:N-subtle}
In higher dimensions the Lusin $(N)$-property is subtle even for mappings in $W^{1,n}_{\rm loc}$:
additional structural hypotheses (e.g.\ openness/discreteness, finite distortion, or suitable modular
inequalities) are often needed to guarantee the validity of condition $(N).$ For background and further results we refer to
\cite[Ch.~4]{HK14} and \cite{KMZ12} and to the references therein.
\end{remark}


\subsection{Jacobian control and area / change-of-variables formulas}\label{subsec:area} 

By Theorem~\ref{thm:ACLp}, every continuous weakly Lipschitz mapping
whose scalar coordinate functions, restricted to almost every
coordinate-parallel line, have Lusin's
$(N)$-property belongs to $W^{1,n}_{\rm loc}(D)$ and hence is
approximately differentiable a.e.
In particular, its Jacobian determinant $J_f$ exists a.e. At points of differentiability one has $\|Df(x)\|_{\rm op}\le L(x,f),$ hence
\begin{equation*}
|J_f(x)|\le \|Df(x)\|_{\rm op}^n \le L(x,f)^n \le \Phi'(x)
\qquad\text{for a.e. }x\in D,
\end{equation*}
where $\|Df(x)\|_{\rm op}$ denotes the operator norm. In particular, the same estimate holds with
$\|Df(x)\|_{\rm op}$ replaced by any matrix norm dominating it (e.g.\ the Hilbert--Schmidt norm).
Consequently, $J_f\in L^1_{\rm loc}(D).$

For a measurable set $A\subset D$ and $y\in\mathbb R^n$, recall that the quantity
\begin{equation*}
N(f,A,y):=\#\{x\in A:\ f(x)=y\}\in\{0,1,2,\dots,\infty\}
\end{equation*}
denotes the multiplicity function.

We record the following standard area (change-of-variables) formula with multiplicity,
valid for Sobolev mappings under Lusin's $(N)$-property; see, e.g., \cite[Ch.~6]{HK14}
and \cite[Ch.~12]{IM01}.

\medskip
\begin{theorem}\label{thm:area-formula} 
{\em Let $f:D\to\mathbb R^n$ be a continuous weakly Lipschitz mapping whose
scalar coordinate functions, restricted to almost every coordinate-parallel line,
have Lusin's $(N)$-property. Assume that $f$ has Lusin's $(N)$-property
in $D.$ Then for every
measurable set $A\subset D$ and every nonnegative measurable function
$\eta:\mathbb R^n\to [0,\infty]$ we have
\begin{equation}\label{eq:area-formula}
\int_A \eta\bigl(f(x)\bigr)\,|J_f(x)|\,dm(x)
  \;=\;
\int_{\mathbb R^n} \eta(y)\,N(f,A,y)\,dm(y).
\end{equation}
In particular, taking $\eta\equiv 1$ yields
\begin{equation*}
\int_A |J_f(x)|\,dm(x)=\int_{\mathbb R^n} N(f,A,y)\,dm(y).
\end{equation*}}
\end{theorem}

\begin{proof}
By Theorem~\ref{thm:ACLp}, $f\in W^{1,n}_{\rm loc}(D,\mathbb R^n).$
The asserted formula is then the standard area formula for Sobolev
mappings satisfying Lusin's $(N)$-property; see \cite{HK14}.
\end{proof}

\medskip
\begin{remark}\label{rem:cov}
If $f$ is injective on $A,$ then $N(f,A,y)=\mathbf 1_{f(A)}(y)$ and \eqref{eq:area-formula}
reduces to the usual change-of-variables formula
\begin{equation*}
\int_A \eta\bigl(f(x)\bigr)\,|J_f(x)|\,dm(x)=\int_{f(A)} \eta(y)\,dm(y).
\end{equation*}
More generally, \eqref{eq:area-formula} gives the change-of-variables formula \emph{with multiplicity}.
\end{remark}

\medskip
\begin{remark}\label{rem:needN}
The role of Lusin's $(N)$-property in Theorem~\ref{thm:area-formula} is not merely technical. Weak Lipschitzness yields $J_f\in L^1_{\rm loc}(D),$ but this alone does not guarantee that null sets have null images. Indeed, if the area formula were valid for every measurable set $A\subset D,$ then applying it to a set $E\subset D$ with $m E=0$ and taking $\eta\equiv 1$ would give
\begin{equation*}
\int_{\mathbb R^n} N(f,E,y)\, dm(y)=0.
\end{equation*}
Since $N(f,E,y)\ge 1$ for $y\in f(E),$ this would imply $m f(E)=0.$ Thus the full area formula itself forces the Lusin $(N)$-property. In this sense, Corollary~\ref{cor:Nproperty} provides the additional null-set control needed to pass from Jacobian integrability to the change-of-variables formula.
\end{remark}


\subsection{Regularity upgrades under stronger integrability of $\Phi'$}\label{subsec:holder} 

For Lipschitz mappings satisfying \eqref{eq:limLip}, the $\beta$-dimensional Hausdorff measure of the image
is controlled by that of the preimage,
\begin{equation}\label{eq:Lip-Hbeta}
\mathcal H^\beta f(A) \le C\,\mathcal H^\beta A,
\end{equation}
see \cite{MV90}. It is necessary to emphasize that such an estimate fails in general for finitely
Lipschitz mappings, and this failure may occur already within the class of weakly Lipschitz
mappings.

\begin{example}\label{ex:MV-weakLip}
Consider the planar mapping \cite[3.3]{MV90}
\begin{equation*}
f(x_1,x_2)=(x_1,g(x_1)),\qquad (x_1,x_2)\in (0,1)^2,
\end{equation*}
where $g:[0,1]\to[0,1]$ is the piecewise linear function defined for each $n=1,2,\dots$ on
the segment $\bigl[\frac1{n+1},\frac1n\bigr]$ by
\begin{equation*}
g(t)=
\begin{cases}
-2n\left(t-\dfrac{2n+1}{2n(n+1)}\right),
& \dfrac{1}{n+1}\le t \le \dfrac{2n+1}{2n(n+1)},\\[2.0ex]
2(n+1)\left(t-\dfrac{2n+1}{2n(n+1)}\right),
& \dfrac{2n+1}{2n(n+1)}\le t \le  \dfrac{1}{n},
\end{cases}
\qquad n=1,2,\ldots,
\end{equation*}
and set $g(0)=0.$ Note that the breakpoint
$c_n=\dfrac{2n+1}{2n(n+1)}$ is the midpoint of every segment $\bigl[\frac1{n+1},\frac1n\bigr].$

Then $f$ is continuous. Since $g$ is piecewise linear, $g'$ is locally
square integrable in $(0,1),$ and hence
\begin{equation*}
L(x,f)\le \sqrt{1+|g'(x_1)|^2}
\quad\text{for a.e. }x=(x_1,x_2)\in D.
\end{equation*}
Consequently $f$ is weakly Lipschitz in $D,$ for instance with
\begin{equation*}
\Phi(A)=\int_A \bigl(1+|g'(x_1)|^2\bigr)\,dm(x),
\qquad A\subset D \text{ open}.
\end{equation*}
However, for the line segment
\begin{equation*}
E=(0,1)\times\{1/2\}\subset D
\end{equation*}
we have
\begin{equation*}
\mathcal H^1 E=1,
\end{equation*}
whereas
\begin{equation*}
f(E)=\{(t,g(t)):0<t<1\}
\end{equation*}
has infinite one-dimensional Hausdorff measure. Hence an estimate of the form
\begin{equation*}
\mathcal H^\beta f(E) \le C\mathcal H^\beta E
\end{equation*}
cannot hold for weakly Lipschitz mappings in general, already for $\beta=1.$
\end{example}

\medskip

A natural positive substitute for \eqref{eq:Lip-Hbeta} follows from improved integrability of $\Phi'.$

\begin{proposition}\label{prop:holder} 
{\em Let $f:D\to\mathbb R^n$ be a continuous weakly Lipschitz mapping whose
scalar coordinate functions, restricted to almost every
coordinate-parallel line, have Lusin's
$(N)$-property. Assume that $\Phi'\in L^p_{\rm loc}(D)$ for some $p>1.$
Then $L(\cdot,f)\in L^{np}_{\rm loc}(D)$ and
$f\in W^{1,np}_{\rm loc}(D).$
In particular, by the Sobolev--Morrey embedding, $f$ has a locally H\"older continuous
representative with exponent
\begin{equation}\label{eq:alpha}
\alpha=1-\frac1p.
\end{equation}}
\end{proposition}

\begin{proof}
Since $L(x,f)^n\le \Phi'(x)$ a.e., the assumption
$\Phi'\in L^p_{\rm loc}(D)$ gives
$L(\cdot,f)\in L^{np}_{\rm loc}(D).$ By Theorem~\ref{thm:ACLp},
$f\in {\rm ACL}^n,$ and at a.e. point
\begin{equation*}
\left|\frac{\partial f}{\partial x_j}\right|\le L(x,f).
\end{equation*}
Hence all first partial derivatives belong to $L^{np}_{\rm loc},$ so
$f\in W^{1,np}_{\rm loc}(D).$ Since $np>n,$ the Sobolev--Morrey
embedding gives the stated local H\"older continuity with exponent
$\alpha=1-n/(np).$
\end{proof}

\begin{corollary}\label{cor:Hdist} 
{\em Under the assumptions of Proposition~\ref{prop:holder}, for every compact set $K\Subset D$
there exists a constant $C=C(K)$ such that
\begin{equation*}
|f(x)-f(y)|\le C|x-y|^\alpha,\qquad x,y\in K,
\end{equation*}
where $\alpha$ is given by \eqref{eq:alpha}. Consequently, for each $s\ge 0$ and each set $E\subset K,$
\begin{equation}\label{eq:Hdist}
\mathcal H^s f(E) \le C_s\,\mathcal H^{\alpha s}E.
\end{equation}}
\end{corollary}

\begin{proof}
For $s=0$ the assertion is immediate, since
$\mathcal H^0 f(E) \le \mathcal H^0 E.$ Assume now that $s>0.$
By Proposition~\ref{prop:holder} and the Sobolev--Morrey embedding,
$f$ is locally $\alpha$-H\"older continuous, where
$\alpha=1-1/p.$ Since the original mapping $f$ is continuous, it
coincides everywhere with its H\"older continuous representative.

Let $K\Subset D.$ Cover $K$ by finitely many balls $B_i$ such that
$2B_i\Subset D.$ For each $i$ there exists a constant $C_i$ such that
\begin{equation*}
|f(x)-f(y)|\le C_i |x-y|^\alpha, \qquad x,y\in B_i .
\end{equation*}
Let $\delta>0$ be a Lebesgue number of this finite cover of $K,$ and put
$C_0=\max_i C_i$ and $M=\operatorname{diam} f(K).$ If $x,y\in K$ and
$|x-y|<\delta,$ then $x,y\in B_i$ for some $i,$ and hence
\begin{equation*}
|f(x)-f(y)|\le C_0 |x-y|^\alpha .
\end{equation*}
If $|x-y|\ge\delta,$ then
\begin{equation*}
|f(x)-f(y)|\le M\le M\delta^{-\alpha}|x-y|^\alpha .
\end{equation*}
Thus the asserted H\"older estimate holds on $K$ with
$C=\max\{C_0,M\delta^{-\alpha}\}.$

It remains to prove the Hausdorff measure estimate. Let $E\subset K,$ and
let $E\subset \bigcup_i U_i$ be a cover with $\operatorname{diam} U_i<\rho.$
Then
\begin{equation*}
f(E)\subset \bigcup_i f(E\cap U_i),
\end{equation*}
and the H\"older estimate on $K$ gives
\begin{equation*}
\operatorname{diam} f(E\cap U_i)
\le C\,(\operatorname{diam} U_i)^\alpha .
\end{equation*}
Therefore, raising the last inequality to the power $s,$ taking
the infimum over all such covers, and then letting $\rho\to 0,$ we get
\begin{equation*}
\mathcal H^s f(E) \le C_s\,\mathcal H^{\alpha s} E.
\end{equation*}
\end{proof}

\begin{remark}\label{rem:sharp}
The exponent $\alpha=1-\frac1p$ in Proposition~\ref{prop:holder} is consistent with the model
one-dimensional map $t\mapsto \operatorname{sgn}(t)|t|^\alpha$, for which the scaling of coverings
shows that the exponent $\alpha$ in \eqref{eq:Hdist} cannot, in general, be replaced by any larger one.
\end{remark}

\begin{remark}\label{rem:Lip-endpoint}
If $\Phi'\in L^\infty_{\rm loc}(D),$ then $L(\cdot,f)\in L^\infty_{\rm loc}(D)$ and $f$ has a locally
Lipschitz representative. In this case, estimates of the form \eqref{eq:Lip-Hbeta} hold on compact
subsets of $D$ for every $\beta\in[0,n].$
\end{remark}


\subsection{Inverse mapping questions for weakly Lipschitz homeomorphisms}\label{subsec:inverse} 

For homeomorphisms it is natural to investigate when $f^{-1}$ is ${\rm ACL}$ (or belongs to a Sobolev
class), and when the Lusin $(N^{-1})$-property holds. In addition to the upper control
\eqref{eq:limset}, a convenient complementary hypothesis is a lower control of the pointwise quantity
$l(x,f)$ defined in \eqref{eq:lowerl}.

\medskip

\noindent\textbf{Two-sided metric control.}
Let $f:D\to D'$ be a homeomorphism. Assume that $f$ is weakly Lipschitz in $D,$ i.e.
\begin{equation*}
L(x,f)\le [\Phi'(x)]^{1/n}\quad\text{for a.e. }x\in D
\end{equation*}
for some $\Phi\in\mathcal F.$
Assume further that there exists a set function $\Upsilon\in\mathcal F$ such that
\begin{equation}\label{eq:lower-control}
l(x,f)\ge [\Upsilon'(x)]^{-1/n}\quad\text{for a.e. }x\in D.
\end{equation}
Equivalently,
\begin{equation*}
\limsup_{y\to x}\frac{|y-x|}{|f(y)-f(x)|}\le [\Upsilon'(x)]^{1/n}\quad\text{for a.e. }x\in D.
\end{equation*}

\begin{proposition}\label{prop:inv-weakLip} 
{\em Let $f:D\to D'$ be a homeomorphism satisfying \eqref{eq:lower-control}, and assume in addition
that $f$ has Lusin's $(N)$-property.
Define
\begin{equation*}
w(y):=\Upsilon'\bigl(f^{-1}(y)\bigr),\qquad y\in D'.
\end{equation*}
Assume that $w\in L^1_{\rm loc}(D')$ and set
\begin{equation*}
\widetilde\Phi(A):=\int_A w(y)\,dm(y),\qquad A\subset D'\ \text{open}.
\end{equation*}
Then $\widetilde\Phi\in\mathcal F$ and $f^{-1}$ is weakly Lipschitz in
$D',$ with
\begin{equation}\label{eq:inv-weakLip}
L(y,f^{-1})\le [\widetilde\Phi'(y)]^{1/n}
\quad\text{for a.e. }y\in D'.
\end{equation}
If, in addition, the scalar coordinate functions of $f^{-1},$ restricted to almost
every coordinate-parallel line in $D',$ have Lusin's $(N)$-property, then
$f^{-1}$ is ${\rm ACL}$ and $f^{-1}\in {\rm ACL}^n.$}
\end{proposition}

\begin{proof}
For $y=f(x)$ we have, by the definition of $L(\cdot,f^{-1})$ and $l(\cdot,f),$
\begin{equation*}
L\bigl(f(x),f^{-1}\bigr)
=\limsup_{z\to f(x)}\frac{|f^{-1}(z)-x|}{|z-f(x)|}
=\limsup_{u\to x}\frac{|u-x|}{|f(u)-f(x)|}
=\frac{1}{l(x,f)}.
\end{equation*}
Thus, \eqref{eq:lower-control} implies
\begin{equation*}
L\bigl(f(x),f^{-1}\bigr)\le [\Upsilon'(x)]^{1/n}\quad\text{for a.e. }x\in D.
\end{equation*}
Let $E\subset D$ be the exceptional set where the last inequality may fail.
Then $m E=0.$ Since $f$ has Lusin's $(N)$-property, we have
$m f(E) =0.$ Hence, for a.e. $y\in D',$ writing $x=f^{-1}(y),$ we obtain
\begin{equation*}
L(y,f^{-1})\le [\Upsilon'(f^{-1}(y))]^{1/n}=[w(y)]^{1/n}.
\end{equation*}
Since $\widetilde\Phi(A)=\int_A w(y)\,dm(y),$ we have
$\widetilde\Phi'(y)=w(y)$ for a.e. $y\in D'.$ Therefore
\begin{equation*}
L(y,f^{-1})\le [w(y)]^{1/n}
=[\widetilde\Phi'(y)]^{1/n}\quad\text{a.e. in }D',
\end{equation*}
which is \eqref{eq:inv-weakLip}. If, in addition, the scalar coordinate functions of $f^{-1},$ restricted to almost
every coordinate-parallel line in $D',$ have Lusin's $(N)$-property, then Theorem~\ref{thm:ACLp}
applied to $f^{-1}$ gives $f^{-1}\in {\rm ACL}^n.$
\end{proof}

\begin{corollary}\label{cor:Ninv} 
{\em In addition to the assumptions of Proposition~\ref{prop:inv-weakLip}, assume that $f^{-1}$
satisfies a Rado--Reichelderfer type diameter estimate
\begin{equation*}
\bigl(\operatorname{diam} f^{-1}(B)\bigr)^n
\le C\,\widetilde\Psi(2B)
\end{equation*}
for every ball $B\subset D'$ such that $\overline{2B}\subset D',$ where
$\widetilde\Psi\in\mathcal F$ is defined on open subsets of $D'.$
Assume also that $\widetilde\Psi$ satisfies the outer absolute-continuity
condition of Corollary~\ref{cor:Nproperty}. Then $f^{-1}$ has Lusin's $(N)$-property
in $D',$ and therefore $f$ has Lusin's $(N^{-1})$-property.}
\end{corollary}

\begin{proof}
Apply Corollary~\ref{cor:Nproperty} to the mapping $g=f^{-1}:D'\to D.$
Since $f$ is a homeomorphism, $g$ is continuous. The assumed estimate
\begin{equation*}
\bigl(\operatorname{diam} f^{-1}(B)\bigr)^n
\le C\,\widetilde\Psi(2B)
\end{equation*}
is precisely the Rado--Reichelderfer diameter estimate required in
Proposition~\ref{prop:image-measure}, applied in the domain $D'$ to $g$
with the set function $\widetilde\Psi$ (absorbing the constant $C$ into
the dimensional constant). The assumed outer absolute-continuity condition
for $\widetilde\Psi$ is exactly the condition used in
Corollary~\ref{cor:Nproperty}. Hence $g=f^{-1}$ has Lusin's
$(N)$-property in $D'.$

Now let $E\subset D'=f(D)$ and suppose that $mE=0.$ Since $f^{-1}$
has Lusin's $(N)$-property, we have
$m f^{-1}(E)=0.$
This is precisely Lusin's $(N^{-1})$-property for $f.$
\end{proof}

\begin{remark}\label{rem:Ninv-need}
The assumption \eqref{eq:lower-control} (or any comparable two-sided metric control) is essential:
without extra hypotheses, the $(N^{-1})$-property may fail even for smooth homeomorphisms
(see, e.g., the example in \cite{KMR19}).
\end{remark}


\subsection{Related remarks and examples}\label{subsec:related}
Here we briefly discuss some other results and further perspectives relating to the class of weakly Lipschitz mappings.


\subsubsection{Upper gradient viewpoint.}
For mappings satisfying the hypotheses of Theorem~\ref{thm:ACLp},
we have $f\in W^{1,n}_{\rm loc}(D).$ Hence the operator norm $\|Df\|_{\rm op}$
is an $n$-weak upper gradient of $f.$ In view of \eqref{eq:limset} we have $\|Df\|_{\rm op}\le \Phi'^{\,1/n}$
a.e., hence $\Phi'^{\,1/n}$ is also an admissible $n$-weak upper gradient in the sense of the
upper-gradient approach to Sobolev mappings; see \cite{HKST15}.


\subsubsection{Remarks on modular / capacity inequalities.}
The idea of controlling a mapping by a measurable density also appears,
in a global form, in the theory of Poletsky-type modulus inequalities and
capacitary distortion. Gehring \cite{Geh71} studied mappings for which the
$p$-capacity of an image ring is bounded by a constant multiple of the
corresponding capacity in the preimage, $p\ne n,$ and showed that such
bounds lead to sharp Lipschitz-type conclusions.

Cristea \cite{Cris14,Cris21} developed related generalized modular
inequalities for local homeomorphisms and for open light mappings in metric
measure spaces. These results are formulated in terms of moduli of curve
families or capacities of condensers, whereas our approach is local; however,
both frameworks use a measurable majorant as a quantitative device for
controlling the mapping's geometric distortion.


\subsubsection{Non-ACL examples and exceptional sets of $\sigma$-finite $\mathcal H^{n-1}$.}
Koskela and Rogovin constructed examples showing that strong metric control
(e.g.\ in terms of linear dilatation) does not by itself imply absolute
continuity on lines; see \cite{KR05}. In such results it is natural to
describe the failure of ${\rm ACL}$ (or of related regularity properties)
via exceptional sets of $\sigma$-finite $(n-1)$-dimensional Hausdorff
measure. This size is geometrically meaningful: in the definition of
${\rm ACL}$ (Section~\ref{subsec:ACL}), for each coordinate direction
$x_j$ the possible failure of absolute continuity is detected on a set
$A_j\subset T_jQ$ of \emph{line parameters}, where
$T_jQ\subset\mathbb R^{n-1}$ is the orthogonal projection of the cube.
Thus the relevant parameter space is $(n-1)$-dimensional. In particular,
codimension-one exceptional sets may interact in an essential way with
families of curves, moduli, and capacities, and in many borderline regularity
statements one cannot, in general, replace $\sigma$-finite
$\mathcal H^{n-1}$ by a higher-codimension exceptional class. This
codimension-one viewpoint is also relevant in refined geometric
characterizations of weak $p$-quasiconformality; see \cite{SU25}.


\subsubsection{Non-homeomorphic mappings and $U$-points.}
Our definition of weak Lipschitzness does not require injectivity. In the non-homeomorphic
setting it becomes natural to compare weakly Lipschitz mappings with the classical theory of
continuous mappings and monogeneity developed by Trokhimchuk \cite{Tro64}. In particular, his
notion of $U$-points provides a language to describe the local structure
of preimages and the ``branching'' behaviour of continuous mappings. Such ideas may be useful when
studying to what extent Lusin-type properties, openness/discreteness, or change-of-variables
formulas can be obtained for weakly Lipschitz mappings without additional topological assumptions.


\Needspace{5\baselineskip}

\bigskip
\noindent
{\bf Declarations.}

Conflict of Interest: None.

Ethical Approval: Not applicable.

Data Availability: No datasets were generated or analyzed during this study.


\bigskip

\noindent{\textbf{Acknowledgement.}}\ \
I am grateful to Pekka Koskela for his helpful comments on this manuscript and for many enjoyable meetings and conversations over the years.



\end{document}